\documentclass[12pt]{article}

\usepackage{amstext}
\usepackage{amsthm}
\usepackage{amsmath}
\usepackage{amssymb}
\usepackage{latexsym}
\usepackage{amsfonts}
\usepackage{color}
\usepackage{graphicx}
\usepackage{yfonts}

\makeatletter
\renewcommand*\env@matrix[1][*\c@MaxMatrixCols c]{%
  \hskip -\arraycolsep
  \let\@ifnextchar\new@ifnextchar
  \array{#1}}
\makeatother

\AtBeginDocument{%
 \def\MR#1{} 
}

\usepackage[pagebackref,hypertexnames=false, colorlinks, citecolor=black, linkcolor=black, urlcolor=red]{hyperref}
\usepackage[backrefs]{amsrefs}
\newcommand\blue{\color{black}}
\newcommand\black{\color{black}}

\usepackage{latexsym}
\usepackage{amssymb}
\usepackage{euscript}

\let\cal=\mathcal      

\def\mcc{M\raise.5ex\hbox{c}C}
\def\mccarthy{M\raise.5ex\hbox{c}Carthy}

\def\eg{{\it e.g. }}
\def\ie{{\it i.e. }}

\def\M{{\cal M}}

\def\z{\zeta}

\def\vare{\varepsilon}

\let\i=\infty

\def\la{\langle}
\def\ra{\rangle}
\def\={\ = \ }

\def\E{E_\l}

\def\C{\mathbb C}

\def\T{\mathbb T}
\def\D{\mathbb D}

\def\be{\setcounter{equation}{\value{theorem}} \begin{equation}}
\def\ee{\end{equation} \addtocounter{theorem}{1}}
\def\beq{\begin{eqnarray*}}
\def\eeq{\end{eqnarray*}}
\def\se{\setcounter{equation}{\value{theorem}}} 
\def\att{\addtocounter{theorem}{1}}
\def\vs{\vskip 5pt}

\def\bp{{\sc Proof: }}
\def\ep{{}{\hfill $\Box$} \vskip 5pt \par}

\def\oec{{}{\hspace*{\fill} $\lhd$} \vskip 5pt \par}
\def\bl{\begin{lemma}}
\def\el{\end{lemma}}
\def\bt{\begin{theorem}}
\def\et{\end{theorem}}
\def\bprop{\begin{prop}}
\def\eprop{\end{prop}}
\def\bd{\begin{definition}}
\def\ed{\end{definition}}
\def\br{\begin{remark}}
\def\er{\end{remark}}
\def\bexer{\begin{exercise}}
\def\eexer{\end{exercise}}

\newtheorem{theorem}{Theorem}[section]
\newtheorem{prop}[theorem]{Proposition}
\newtheorem{lemma}[theorem]{Lemma}
\newtheorem{cor}[theorem]{Corollary}

\newtheorem{definition}[theorem]{Definition}
\newtheorem{examp}[theorem]{Example}

\newtheorem{defin}[theorem]{Definition}

\renewcommand\Re{\mathrm{Re\, }}

\renewcommand\O{\Omega}
\def\gdel{{\mathcal G}_\delta}

\def\L{{\mathcal L}}

\def\vare{\varepsilon}

\def\M{{\mathbb M}}

\def\U{\mathcal U}

\def\mn{\M_n}

\def\mmd{\M_m^d}
\def\mnd{\mn^d}

\def\md{{\mathbb M}^{[d]}}

\def\d{\delta}
\def\ph{\varphi}

\newcommand{\tensor}[2]{\text{ }{\begin{smallmatrix} #1 \\ \otimes\\ #2\end{smallmatrix}}\text{  }}

\newcommand{\tensthree}[3]{\text{ }{\begin{smallmatrix} #1 \\ \otimes\\ #2 \\ \otimes \\ #3 \end{smallmatrix}}\text{  }}

\newcommand\yt{Y_T}

\renewcommand\E{{\mathcal E}}
\newcommand\idel{{\mathcal I}_\delta}
\newcommand\bdel{{\mathcal B}_\delta}
\newcommand\sgd{{\mathcal S}(\gdel)}
\newcommand\tp{T}

\numberwithin{equation}{section}
\title{A non-commutative Julia Inequality}
\author{John E. M\raise.5ex\hbox{c}Carthy
\thanks{Partially supported by National Science Foundation Grant  
DMS 1565243}
\and
James E. Pascoe
\thanks{Partially supported by National Science Foundation Fellowship  
DMS 1606260}
}

\begin{document}

\bibliographystyle{plain}
\maketitle

\begin{abstract}
We prove a Julia inequality for bounded non-commutative functions on polynomial polyhedra.
We use this to deduce a Julia inequality for holomorphic functions on classical domains in $\C^d$.
We look at differentiability at a boundary point for functions that have a certain regularity there.
\end{abstract}

\section{Introduction}

The classical Julia inequality asserts that if a holomorphic function $\phi$ maps the unit disk $\D$ to itself,
and if at some boundary point $\tau \in \partial \D$ one has
\be\label{eqa1x}
\liminf_{z\to\tau} \frac{1-|\ph(z)|^2}{1-|z|^2} \= \alpha \ < \infty,
\ee
 then there exists $\omega \in \partial \D$ such that
\be
\label{eqa2}
\frac{|\ph(z) - \omega |^2}
{1-|\ph(z)|^2}
\ \leq \
\alpha \ \frac {|z - \tau|^2} 
{1-|z|^2} .
\ee
The inequality was proved, with an extra regularity hypothesis on $\ph$, by G. Julia
in \cite{ju20}, and in the form stated by C. Carath\'eodory in \cite{car29}.
D. Sarason found a proof using model theory \cite[Chap VI]{sar94}.

Generalizations of Julia's inequality have been found for functions on the ball by M. Herv\'e \cite{her63}, W. Rudin \cite[Sec. 8.5]{rud76} and M. Jury \cite{ju07},
and on the polydisk by K.~Wlodarczyk \cite{wlo87},  F.~Jafari  \cite{jaf93}  and M.~Abate \cite{ab98}.
In the case of the bidisk,  a detailed analysis of points for which the analogue of the Julia quotient \eqref{eqa1x} 
remains bounded (these are called B-points)  has been carried out in \cite{amy12a, amy11b,aty12,bk13}.

It is the purpose of this note to extend Julia's inequality, and the study of B-points, to 
non-commutative functions (which we shall define in Subsection \ref{seca1} below)
that are bounded on polynomial polyhedra. Our methods rely 
on the model-theoretic
ideas of \cite{amy12a}. 

Our results are of interest even
in the commutative case,
because they provide a unified approach to 
proving boundary versions of the Schwarz-Pick lemma in the Schur-Agler class of various domains, such as
the polydisk, or the multipliers of the Drury-Arveson space. The methods also show that at B-points where the function
is not analytic, it does have directional derivatives in all directions pointing into the set, and the derivative is 
a holomorphic (but not necessarily linear) function of the direction.
This is explained in Section~\ref{secc} below.

\subsection{ Non-commutative Functions}
\label{seca1}

Non-commutative function theory, which originated in the work of J.L. Taylor
\cite{tay72,tay73}, has recently started to flourish 
 \cite{po06,po11,bgm06,akv06,hkm11b,hm12,ms13,chmn13,mt15,pas14,ptd17,bal15,amco16}.
The foundations are developed in the book \cite{kvv14}.
 
The idea is to study functions of non-commuting variables that are generalized  non-commuting polynomials
in an analogous fashion to thinking of a holomorphic function as a generalization of a polynomial.
Our domains are domains of $d$-tuples of matrices, but they don't reside in just one dimension.
We let $\mn$ denote the $n$-by-$n$ complex matrices, and let
\[
\md \= \cup_{n=1}^\i \mnd .
\]
We shall call a function $\phi$  defined on a subset of $\md$ {\em graded} if, whenever 
$x \in \mnd$, then $\phi(x) \in \mn$. If $x \in \mnd$ and $y \in \mmd$, we shall let $x \oplus y$ denote
the $d$-tuple in $\M_{n+m}$ obtained by direct summing each component.
If $x \in \mnd$ and $s$ is an invertible matrix in $\mn$, then
$s^{-1} x s$  denotes the $d$-tuple  $(s^{-1} x^1 s, \dots, s^{-1} x^d s) $.

\bd
\label{defa1}
 An nc-function $\phi$ on a set $\Omega \subseteq \md$ is a graded function that respects direct sums and joint similiarities, \ie 
\beq
\phi(x \oplus y) &\= &  \phi(x) \oplus \phi(y)\\
\phi(s^{-1} x s) := \phi (s^{-1} x^1 s, \dots, s^{-1} x^d s) &= & s^{-1}\phi(x) s, 
\eeq
where the equations are only required to hold when the arguments on both sides are in $\Omega$,
and in the second one if $x \in \mnd$, then $s$ is invertible in $\mn$.
\ed

(We  use superscripts for components, since we shall have many sequences indexed by subscripts.)
Notice that every non-commutative polynomial is an nc-function on $\md$. 
A particularly nice class of domains on which to study nc-functions are polynomial polyhedra.
These are defined in terms of a matrix $\delta$, each of whose entries is a non-commutative polynomial in $d$ variables.
Then $\gdel$ is defined as
\be
\label{eqab01}
\gdel \ := \ \cup_{n=1}^\infty \{ x \in \mnd \ : \ \| \delta(x) \| < 1 \} ,
\ee
\blue
where $\d(x)$,  a matrix of matrices, is given the operator norm.
\black
A primary example is
 the $d$-dimensional noncommutative polydisk, which is the set
\be
\label{eqab2}
\{ x \in \md : \| x\| := \max_{1 \leq r \leq d} \| x^r \| < 1 \} .
\ee
For this set, we can take $\delta$ to be the diagonal $d$-by-$d$ matrix with the coordinate
functions on the diagonal.
Another well-studied example (see \eg \cite{po91}) is the non-commutative ball, 
which  we shall take to be the column contractions\footnote{It is more common to consider the row contractions,
but we choose column contractions so that what we call the distinguished boundary will be non-empty.
It is easy to pass between these two sets, since the column contractions are just the adjoints of the row contractions.},
\be
\label{eqab3}
\{ x \in \md \ : \ \sum_{r=1}^d x^{r*} x^{r} \leq I \},
\ee
which is obtained by letting
\[
\delta(x) \=
\begin{pmatrix}
  x^1 \\ x^2 \\ \vdots \\ x^d 
  \end{pmatrix} .
\]


\subsection{Principal Results}
Let $\gdel$ be defined by \eqref{eqab01}.
We shall let $\bdel$ denote the topological boundary of $\gdel$.
If $x \in \bdel$, then $\| \d(x) \| = 1$, but the converse need not hold ---
\eg with $d=1$, take $\d(x) = \begin{pmatrix}x-1 & 0 \\0 & x+1
\end{pmatrix}$. Then $\gdel$ is empty, but $\| \d(0) \| = 1$.

 The distinguished boundary of $\gdel$, which we shall
denote $\idel$, is 
\be
\label{eqc1}
\idel \ := \ \{ x \in \bdel \ : \ \d(x) \ {\rm is\ an\ isometry} \} .
\ee
The reader should keep in mind the example of the non-commutative polydisk \eqref{eqab2}, in which case
$\idel$ is $\U^{[d]}$, the set of $d$-tuples of unitary matrices in $\md$,
and $\bdel$ is the set of contractive $d$-tuples such that at least one  has norm equal to one.
For the column-ball \eqref{eqab3}, the distinguished boundary will agree with $\bdel$ when $n=1$, but
will be smaller when $n >1$.

\bd
The Schur class of $\gdel$, denoted $\sgd$, is the set of nc functions $\phi$  on $\gdel$ 
such that $\| \phi(x) \| < 1 \ \forall \ x \in \gdel$.
\ed

A B-point for $\phi$ is a point in the  boundary where a certain regularity occurs.
\begin{defin}  
Let $\phi \in \sgd$, and let $T \in \bdel$.
Then $T$ is a B-point for $\phi$ if
\be
\label{eqb3}
\lim_{j \to \i} \frac{ \| I - \phi(Z_j)^* \phi(Z_j) \|}{1 - \|\d(Z_j) \|^2} 
 \ < \ \i 
\ee
 for some sequence $Z_j$ converging to $T$.
\end{defin}

Here is our non-commutative Julia inequality, which says that at a B-point, one has
a boundary version of the Schwarz-Pick inequality, akin to \eqref{eqa2}.
\vs

{\bf Theorem~\ref{thmju}}
{\em
Suppose $\phi \in \sgd$ and $T \in \bdel$. If
\[
\liminf_{Z \to T} \frac{ \| I - \phi(Z)^* \phi(Z) \|}{ 1 - \| \d( Z) \|^2} = \alpha ,
\]
then there exists $W \in \U_n$ such that for all $Z$ in $\gdel$ 
\[
\frac{\| \phi(Z) - W \|^2}{\| I - \phi(Z)^* \phi(Z) \|}
\ \leq \ 
\alpha 
\left( \frac{\| I - \d(T)^* \d(Z) \|^2}{ 1- \| \d(Z) \|^2}\right).
\]
}
\vskip 5pt
{\blue
We note that the noncommutative structure
then gives
  that if $T$ is a $B$-point, then so is $U^*TU$  for any unitary matrix $U$. Furthermore, given two $B$-points $S$ and $T,$ then $S \oplus T$ must also be a $B$-point.
}

Our second  main result, Theorem~\ref{thmf2}, gives a characterization 
of when a point in $\idel$ is a B-point, in terms of a realization of a $\d$ nc-model for $\phi$.
We shall defer an exact statement until Section~\ref{secf}.

Our third result, Theorem~\ref{thmg1}, holds under the assumption that there are a lot of inward directions
at $T$. For now, we shall just give a special case.

\bt \label{specialdirder}
Suppose $\gdel$ is either \eqref{eqab2} or \eqref{eqab3}. 
 Suppose $T \in \idel$ is a B-point of $\phi$.
Then
\[
\eta(H) \= \lim_{t \downarrow 0} \frac{\phi(T + t H) - \phi(T)}{t}
\]
exists for all $H$ satisfying $T + tH \in \gdel$ for $t$ small and positive.
Moreover $\eta$ is a holomorphic function of $H$, which is homogeneous of degree 1.
\et
{
\blue
There is compatibility between the directional derivative at a B-point $T$, given by the function $\eta$ in Theorem \ref{specialdirder},  and the directional derivative at a direct sum of $m$ copies of $T$, which we will denote by $T^{(m)}$. In particular, if a scalar point is a B-point, then one obtains B-points at all levels.
\begin{cor} \label{specialdirdercorollary}
Suppose $\gdel$ is either \eqref{eqab2} or \eqref{eqab3}. 
 Suppose $T \in \idel \cap \mathbb{M}_1^d$ is a B-point of $\phi$.
Then
\[
\eta(H) \= \lim_{t \downarrow 0} \frac{\phi(T^{(m)} + t H) - \phi(T^{(m)})}{t}
\]
exists for all $m,$ and for all $H$ satisfying $T^{(m)} + tH \in \gdel$ for $t$ small and positive.
Moreover $\eta$ is a free function of $H$, which is homogeneous of degree $1$.
\end{cor}
Corollary \ref{specialdirdercorollary} follows immediately from
the definition of a free function, and the fact that $T$ is a $d$-tuple of scalars.
  We leave the details to the interested reader, but emphasize  that the noncommutative Julia-Carath\'eodory really captures a regularity that holds not only within each level  $G_\delta \cap \mathbb{M}_n^d,$ but also between levels.
}

In Example \ref{exg1} we consider the function
\[
\phi(Z^1, Z^2) \= \frac 12 (Z^1 + Z^2) + \frac 12 (Z^1 - Z^2) (2 - Z^1 - Z^2)^{-1} (Z^1 - Z^2)
\]
which we show is in the Schur class of \eqref{eqab01}.


\section{ Background material}
\label{secab}

\subsection{The one variable Julia-Carath\'eodory theorem}
\label{secaba}

 The Julia-Carath\'eodory Theorem, due to G. Julia \cite{ju20} in 1920 and C. Carath\'eodory \cite{car29} in 1929,
is the following.
\bt
\label{jc}
Let $\ph : \D \to \D$ be a holomorphic function, and $\tau \in \T$.
The following are equivalent:

 (A) $
\liminf_{z\to\tau} \frac{1-|\ph(z)|^2}{1-|z|^2} < \infty .$

(B) The quotient
$ \frac{1-|\ph(z)|^2}{1-|z|^2} $
has a non-tangential limit as $z$ tends to $\tau$. 

(C) The function $\ph$ has both a non-tangential limit $\omega \in \T$ at $\tau$ and also an angular derivative $\eta \in \C$, that is the difference quotient
\[
\frac{\ph(z) - \omega}{z - \tau}
\]
has a non-tangential limit $\eta$ at $\tau$.

(D) There exist $\omega$ in $\T$ and $\eta$ in $\C$ so that at $\tau$, $\ph(z)$ tends to $\omega$ non-tangentially 
and $\phi'(\z)$ tends to $\eta$ non-tangentially.

 Furthermore, if \eqref{eqa1x}
  holds, then \eqref{eqa2} does.
 \et

On the bidisk, the analogue of  {\em (B)} does not imply {\em (C)}; but it is proved in \cite{amy12a} that
{\em (C)} implies {\em (D)}. Moreover, it is shown that even when $\phi$ does not have a holomorphic differential
pointing into the bidisk, the one sided derivative exists and is holomorphic in the direction.

\subsection{Background on free holomorphic functions}
\label{secabb}

A free holomorphic function is an nc function that is locally bounded with respect
to the topology generated by all the sets $\gdel$, as $\delta$ ranges over
all matrices with entries that are free polynomials.
These functions are studied in \cite{amfree}, and two principal results are obtained.
One is that a bounded function on $\gdel$ is nc if and only if
it is the pointwise limit of a sequence of non-commutative polynomials.
The other is 
 that every bounded nc-function on $\gdel$ has an nc $\delta$-model.
\blue
 Alternative proofs of both these results have been found in \cite{bmv16b} and \cite{amfreeII}.
 \black
Before explaining what this is, we need to slightly expand definition \ref{defa1}.
Let $\E_1$ and $\E_2$ be Hilbert spaces, and let $\L(\E_1, \E_2)$ denote the bounded linear operators
from $\E_1$ to $\E_2$. Following \cite{ptd17}, we shall \emph{write tensor products vertically} to enhance readability and condense realization formulas,  so
$\tensor{A}{B}$ represents the same object as $A \otimes B$.
We shall assume that the domain $\Omega$ of any nc function is closed w.r.t. direct sums.
Also, we shall for notational convenience assume that $\delta$ is  a square $J$-by-$J$ matrix --- we can always
add rows or columns of zeroes to ensure this.
\blue
We shall let $\O_n $ denote $ \Omega \cap \mnd$, and ${\cal I}_n$ denote the invertible matrices in $\mn$.
\black
\bd
An $\L(\E_1,\E_2)$-valued nc function $F$ 
on a set $\Omega \subseteq \md$ is a  function satisfying
\beq
F(x) \ &\in& \ \L (\tensor{ \E_1}{ \C^n},  \tensor{ \E_2}{ \C^n}) \quad \forall \ n, \ \forall  x \in \Omega_n \\
F(x \oplus y) &=&   F(x) \oplus  F(y)  \quad \forall\  m,n \ \forall \ x \in \Omega_n,\  y \in \Omega_m  \\
F(s^{-1} x s) &=&  \tensor{I_{\E_2}}{s^{-1}} F(x)  \tensor{I_{\E_1}}{s} \quad \forall \ n, \
\forall \ s \in {\cal I}_n \ {\rm s.t.\ }  x, s^{-1} x s \in \Omega_n.
\eeq
\ed

\bd
Let $\delta$ be a  $J$-by-$J$ matrix of free polynomials. Let $\phi$ be an nc-function on $\gdel$.
A $\delta$ nc-model for $\phi$ is an $\L(\C, \tensor{\E}{ \C^J})$-valued nc function $u$ on $\Omega$ that satisfies
\be
\label{eqab1}
I_{\C^n} - \phi(y)^* \phi(x) \= u(y)^* \left(
\tensor{I_\E}{
I_{\C^J \otimes \C^n} -  \delta(y)^* \delta(x)}\right) u(x) .
\ee
\ed

\bt
\label{thmmod}
\cite{amfree}
An nc function $\phi$ defined on $\gdel$ is bounded by $1$ in norm if and only if
it has a $\delta$ nc-model.
\et

\blue
To help the reader, let us rewrite Theorem~\ref{thmmod} in the case that $\gdel$ is the non-commutative $d$-polydisk,
given by \ref{eqab2}, and both $\E_1$ and $\E_2$ are one dimensional. Then $J = d$, and  Theorem~\ref{thmmod}
becomes:

\bt
\label{thmmodpoly}
Let $\O$ be the  non-commutative $d$-polydisk, \eqref{eqab2}.
The graded function $\phi$ is in the Schur class of $\O$  if and only if  there are $d$ 
 $\L(\C, \E)$-valued nc functions $u^1, \dots u^d$ so that, for all $n$, for all $x,y \in \O \cap \mnd$, we have
 \begin{eqnarray*}
\lefteqn{ I_{\C^n} - \phi(y)^* \phi(x)  \= }
\\
&
 \begin{pmatrix}
u^1(y) \\
u^2(y) \\
\vdots \\
u^d(y)
\end{pmatrix}^*
 \begin{pmatrix}
 I -  \tensor{I_\E}{(y^1)^* x^1} & 0 & \dots & 0 \\
 0 &  I -  \tensor{I_\E}{(y^2)^* x^2} & \dots & 0 \\
 \vdots & \vdots & \ddots& \vdots \\
 0 & 0 & \cdots &  \tensor{I_\E}{ I - (y^d)^* x^d}
 \end{pmatrix}
\begin{pmatrix}
u^1(x) \\
u^2(x)\\
\vdots \\
u^d(x)
\end{pmatrix}\\
\\
&= \sum_{j=1}^d u^j(y)^* \left(I_{\E \otimes \C^n} - \tensor{I_\E}{(y^j)* x^j} \right) u^j(x).
\end{eqnarray*}
\et

\black

\section{Julia's Inequality and Consequences}

We shall assume for the remainder of the paper that $\d$ is a $J$-by-$J$ matrix of non-commutative polynomials,
that $\phi \in \sgd$, and $u$ is a $\d$ nc-model for $\phi$,
with values in  $\L(\C, \tensor{\E}{ \C^J})$. We shall further assume that $T \in \bdel$ is in $\mnd$.
We shall let ${\mathcal U}_n$ denote the $n$-by-$n$ unitaries.

\blue
If $\phi$ is in $\sgd$, it follows from Theorem~\ref{thmmod} that 
\[
\| I - \phi(x)^* \phi(x) \| \ \geq \ \|u(x)\|^2 (1 - \| \delta(x)\|^2) .
\]
So if \eqref{eqb3} holds, then $\{ \| u(Z_j) \| \}$ is bounded, and the following definition is non-vacuous.
\black

\begin{defin}
\label{define1}
Let $T$ be a  B-point for $\phi$.
 Let $\yt = \yt(u)$ denote the set of all weak-limits of $u(Z_j)$,
where $Z_j$ is a sequence in $\gdel$ that converges to $T$ and has
\be
\label{eqej1}
\frac{ \| I - \phi(Z_j)^* \phi(Z_j) \|}{ 1 - \|\d( Z_j) \|^2} 
\ee
bounded. 
We shall call $\yt$ the cluster set of the model $u$. 
\end{defin}

\bprop
\label{propb2}
Let $T$ be a  B-point for $\phi$. Then there exists $W \in \U_n$ such that for all $v \in \yt$, for
all $Z$ in $\gdel$, we have
\be
\label{eqb4}
I - W^* \phi(Z) \= v^*
 \left(
\tensor{I}{I -  \delta(T)^* \delta(Z)}\right)
 u(Z) .
\ee
Moreover, if 
\be
\label{eqb5}
\lim_{j \to \i} \frac{ \| I - \phi(Z_j)^* \phi(Z_j) \|}{ 1 - \| \d(Z_j) \|^2} 
\ \leq \alpha 
\ee
 holds for some sequence, then there exists $v \in \yt$ with $\| v\|^2 \leq \alpha$.
\eprop
\bp
Suppose \eqref{eqb3} holds and $u(Z_j)$ tends weakly to $v$. Then 
$\| I - \phi(Z_j)^* \phi(Z_j) \| \to 0$, so by passing to a subsequence we can assume that 
$\phi(Z_j)$ tends to some unitary $W$.
 Taking the limit
in
\[
I - \phi(Z_j)^* \phi(Z) \= u(Z_j)^*
 \left( 
\tensor{I}{I -  \delta(Z_j)^* \delta(Z)}\right)
 u(Z) .
\]
we get
\eqref{eqb4}.
To see that $W$ is unique, if another sequence $Z_j'$ tending to $T$ with $u(Z_j')$ weakly convergent 
had $\phi(Z_j') \to W'$, then letting $Z = Z_j'$ in \eqref{eqb4} 
\blue
we get
\[
I - W^* \phi(Z_j') \= v^*
 \left(
\tensor{I}{I -  \delta(T)^* \delta(Z_j')}\right)
 u(Z_j') .
\]
Since $\phi(Z_j)'$ and $\delta(Z_j')$ act on  finite dimensional spaces, we can pass to a subsequence 
so that they converge in norm, and $u(Z_j')$ converges weakly.
\black
In the limit, we get
\[
I - W^* W' \= 0 ,
\]
so $W = W'$.

For the latter part, note
\beq
(1 - \| \d(Z_j)  \|^2 )u(Z_j)^* u(Z_j)
& \ \leq \ &
u(Z_j)^* 
 \left( 
\tensor{I}{I -  \delta(Z_j)^* \delta(Z_j)}\right)
u(Z_j)
\\
&=& I - \phi(Z_j)^* \phi(Z_j) ,
\eeq
so
\[
\| u(Z_j) \|^2 \ \leq \ \frac{ \| 1 - \phi(Z_j)^* \phi(Z_j) \|}{ 1 - \| \d( Z_j) \|^2} .
\]
Taking $v$ to be any weak cluster point of $u(Z_j)$, we get $\|v\|^2 \leq \alpha$.
\ep

If $T$ is a B-point for $\phi$, we shall let $\phi(T)$ denote the matrix $W$ that satisifies \eqref{eqb4}.

Here is the nc Julia inequality.
\bt
\label{thmju}
Suppose  $T \in \bdel$, and 
\[
\liminf_{\gdel \ni Z \to T} \frac{ \| I - \phi(Z)^* \phi(Z) \|}{ 1 - \| \d( Z) \|^2} = \alpha .
\]
Then there exists $W \in \U_n$ such that for all $Z$ in $\gdel$ 
\be
\label{eqb6}
\frac{\| \phi(Z) - W \|^2}{\| I - \phi(Z)^* \phi(Z) \|}
\ \leq \ 
\alpha 
\left( \frac{\| I - \d(T)^* \d(Z) \|^2} { 1 - \| \d(Z)\|^2 }\right).
\ee
\et
\bp
By Proposition \ref{propb2}, we can choose $v$ in $\yt$ with $\| v \|^2 \leq \alpha$.
From \eqref{eqb4} we have
\[
I - W^* \phi(Z) \=  v^* 
 \left(   \tensor{I}{I -  \delta(T)^* \delta(Z)}\right)
 u(Z)  ,
\]
so
\se\att
\begin{eqnarray}
\nonumber
\|  \phi(Z)  - W \|^2 & \ = \ & 
\|  v^* 
 \left(  \tensor{I}{I -  \delta(T)^* \delta(Z)}\right)
 u(Z)  \|^2 
\\
\label{eqb8}
& \leq & 
\alpha \ \|u(Z) \|^2 \| I -  \delta(T)^* \delta(Z)\|^2 .
\end{eqnarray}
Now
\se\att
\begin{eqnarray}
\nonumber
\| I - \phi(Z)^* \phi(Z) \| &\=&
\|   u(Z)^*  \left(  \tensor{I}{I -  \delta(Z)^* \delta(Z)}\right)
 u(Z) \| \\
&\ \geq \ &
  \|  u(Z) \|^2 (1  - \| \d(Z) \|^2)  .
 \label{eqb9}
\end{eqnarray}
Combining \eqref{eqb8} and \eqref{eqb9}, we get
\[
\|  \phi(Z)  - W \|^2 \ \leq \
\alpha   \left( \frac{\| I - \d(T)^* \d(Z) \|^2}{1 - \| \d(Z)\|^2}\right)  \| I - \phi(Z)^* \phi(Z) \| ,
\]
which yields \eqref{eqb6}. \ep


If $T \in \idel$, a non-tangential approach region is a region of the form
\be
\label{eqa1}
\{ Z : \| \d(Z) - \d(T) \| \leq c (1 - \| \d(Z) \|^2) \} .
\ee
A corollary of Julia's lemma is that a function's behavior is controlled non-tangentially at a B-point on the distinguished boundary.
\bprop
\label{propb3}
If $T \in \idel$ is a $B$-point for $\phi$, then
\[
\frac{\| I - \phi(Z)^* \phi(Z) \|}{1 - \| \d(Z)\|^2}
\]
is bounded on all sets that approach $T$ non-tangentially.
%

\eprop
\bp
Suppose $v$ and $W$ are such that
\eqref{eqb4} holds:
\[
I - W^* \phi(Z) \=  v^{*}
 \left( 
\tensor{I}{ I - \delta(T)^* \delta(Z)}\right) u(Z) .
\]
Fix $c$, and let $S$ be the non-tangential approach region
\[
S \= \{ Z : \| \d(T) - \d(Z) \| \leq c (1 - \|\d(Z)\|^2) \}.
\]
By \eqref{eqb6}, we have for $Z$ in $S$ that
\be
\label{eqb11}
\| \phi(Z) - W \|^2 \ \leq \
\alpha  \ 
\left( \frac{\| I - \d(T)^* \d(Z) \|^2}{1 - \| \d(Z)\|^2}\right)  \| I - \phi(Z)^* \phi(Z) \| .
\ee
Now
\se\att
\begin{eqnarray}
\nonumber
 \| I - \phi(Z)^* \phi(Z) \| &\ = \ &
 \| W^* W -  \phi(Z)^* \phi(Z) \| \\
 \nonumber
 & = & \| W^*( W - \phi(Z) ) - (\phi(Z)^* - W^*) \phi(Z) \| \\
 \nonumber
 &\leq & \| W - \phi(Z) \| (\| W\| + \| \phi(Z) \|) \\
 \label{eqbf2}
 &\leq & 2 \| W - \phi(Z) \| .
 \end{eqnarray}
 Squaring and using \eqref{eqb11}, we get
 \be
\label{eqb12}
  \| I - \phi(Z)^* \phi(Z) \|^2 \ \leq \
  4  \alpha 
\left( \frac{\| I - \d(T)^* \d(Z) \|^2}{1 - \| \d(Z)\|^2}\right)  \| I - \phi(Z)^* \phi(Z) \| .
  \ee
Since $T \in \idel$, we have $\d(T)$ is an isometry, so
\[
\| \d(T) - \d(Z) \| \= \| I - \d(T)^* \d(Z) \| .
\]
  Therefore if $Z \in S$, the expression in parantheses on the right-hand side of \eqref{eqb12}
is bounded by
$c^2 (1 - \| \d(Z) \|^2) $, so we get
\[
 \| I - \phi(Z)^* \phi(Z) \| \ \leq \
4 \alpha  c^2  (1 - \|\d( Z)\|^2) ,
\]
as required.
\ep

\section{Models and B-points}
\label{secf}

In this section we shall study how being a B-point is related to properties of the $\d$ nc model.
In the case of the bidisk, these results  are in \cite{amy12a}.

\bprop 
\label{propf1}
Let $T \in 
\idel$,
 and suppose $Z_j$ in $\gdel$ converges to $T$ non-tangentially in the region
\eqref{eqa1}.
The following are equivalent:

(i) $ \| I -  \phi(Z_j)^* \phi(Z_j) \| \leq M  \|I - \d(Z_j)^* \d(Z_j)\|.$

(ii)  $ \| I -  \phi(Z_j)^* \phi(Z_j) \| \leq M' ( 1 - \| \d(Z_j)\|^2).$

(iii) $\| u(Z_j) \|^2 \leq M''$ for some $u$ satisfying \eqref{eqab1}.

(iv) $\| u(Z_j) \|^2 \leq M''$ for every $u$ satisfying \eqref{eqab1}.
\eprop

\bp

$(i) \Rightarrow (iv)$ 
\beq
(1 - \|\d( Z_j)  \|^2 )u(Z_j)^* u(Z_j)
& \ \leq \ &
u(Z_j)^*
 \left( 
\tensor{I}{I - \delta(Z_j)^* \delta(Z_j)}\right)
 u(Z_j) \\
&=& I - \phi(Z_j)^* \phi(Z_j) \\
& \ \leq \  &
M \| I - \d(Z_j)^* \d(Z_j) \| \ I.
\eeq
Therefore
\[
(1 - \| \d(Z_j)  \|^2 )\| u(Z_j) \|^2 \ \leq \
M \| I - \d(Z_j)^* \d(Z_j) \| .
\]
But
\beq
\| I - \d(Z_j)^* \d(Z_j) \| 
&\= & 
\| \d( T)^*\d(T) -\d( Z_j)^*  \d(Z_j) \|  \\
&=&
\| \d(T)^* (\d(T) - \d(Z_j)) + (\d(T)^* - \d(Z_j)^*) \d(Z_j) \|
\\
&\leq &
(1 + \| \d(Z_j) \|)(\| \d(T) - \d(Z_j) \|) \\
&\leq & c (1 + \| \d(Z_j)\|) (1 - \| \d(Z_j)\|^2).
\eeq
Therefore we get
\beq
\| u(Z_j) \|^2 & \ \leq \ & M c (1+\| \d(Z_j) \| ) \\
&\leq & 2Mc .  \qquad \hfill \lhd
\eeq

$(iii) \Rightarrow (i)$ Since
\[
I - \phi(Z_j)^* \phi(Z_j) \= u(Z_j)^*
 \left( 
\tensor{I}{I -  \delta(Z_j)^* \delta(Z_j)}\right)
  u(Z_j), \]
taking norms we get
\[
\| I - \phi(Z_j)^* \phi(Z_j) \| \ \leq \ M'' \| I - \d(Z_j)^* \d(Z_j) \| .
\qquad \hfill \lhd
\]

$(i) \Leftrightarrow (ii)$
In the region \eqref{eqa1}, we have
\beq
1 - \|\d( Z)\|^2 & \ \leq \ &
\| I - \d(Z)^* \d(Z) \| \\
&\leq & 
2c ( 1 - \| \d(Z)\|^2) .
\eeq
\ep

We can now give a different characterization of B-points that are on the distinguished boundary.
\begin{cor}
\label{corf2}
Let $T \in \idel$, and let $u$ be a $\delta$ nc-model for $\phi$. 
Suppose:

(NT)  There exists some sequence in $\gdel$ that approaches $T$ non-tangentially.

Then
The following are equivalent:

(i) $T$ is a B-point of $\phi$.

(ii) $u(Z_j)$ is bounded on some sequence $Z_j$ that approaches $T$ non-tangentially.

(iii) $u(Z)$ is bounded on every set that approaches $T$ non-tangentially.

(iv) $\frac{ \| I - \phi(Z)^* \phi(Z) \|}{ 1 - \| \d(Z) \|^2}$ is bounded on every set
that approaches $T$ non-tangentially.
\end{cor}
\bp
(i) $\Rightarrow$ (iv) by Proposition~\ref{propb3}, and (iii)  $\Rightarrow$  (ii) is trivial.

(iv)  $\Rightarrow$ (iii) and (ii)  $\Rightarrow$ (i) both follow from Proposition~\ref{propf1},
 and the observation that the proof shows that all
the constants $M, M', M''$ are comparable once the aperture of the non-tangential approach region is fixed.
\ep
Remark: Condition (NT) is very mild. It will hold if 
$\Gamma(T)$ (see Def. \ref{defing1} below) is non-empty.


\vskip 5pt
If $u$ is a $\d$ nc model for $\phi$, then by \cite[Cor. 8.2]{amfree}, there is an isometry
(which is called a {\em realization of the model})
\be
\label{eqf9}
\begin{pmatrix}
A&B\\
C&D
\end{pmatrix} \ : \ \C \oplus \tensor{\E}{\C^J} \ \to
\ \C \oplus \tensor{\E}{\C^J}
\ee
so that for $x \in \gdel \cap \mnd$,
\be
\label{eqf1}
\left[
\tensthree{I_\E}{I_{\C^J}}{I_{\C^n}} 
\ - \
\tensor{D}{I_{\C^n}} \ \tensor{I_\E}{\d(x)} \right]
\ u(x) \=
\tensor{C}{I_{\C^n}},
\ee
and
\se
\att
\begin{eqnarray}
\nonumber
\phi(x) &\=& 
\tensor{A}{I_{\C^n}} + \tensor{B}{I_{\C^n}}
\tensor{I_\E}{\d(x)} u(x) \\
\label{eqf10}
&=&
\tensor{A}{I_{\C^n}} + \tensor{B}{I_{\C^n}}
\tensor{I_\E}{\d(x)} \left[ I - \tensor{D}{I_{\C^n}} \tensor{I_\E}{\d(x)} \right]^{-1} 
\tensor{C}{I_{\C^n}}  .
\end{eqnarray}

For $T \in \idel$,  the inward directions for $T$ are those $H$ such that $T + tH$ is inside
$\gdel$ for $t$ small and positive. Formally, if $T \in \mnd$, and $H \in \M_{n}^d$, let
\[
\nabla \delta( \tp) [H] \= \lim_{t \to 0} \frac{1}{t} [ \d( \tp + t H) - \d(\tp) ]
\]
denote the derivative of $\delta$ at $\tp$ in the direction $H$.
 If $A$ is a self-adjoint matrix, we write $A < 0$ to mean $A$ is negative definite.

\begin{defin}
\label{defing1}
Let $T \in \mnd$ be in $\idel$. The inward set of $T$ is the set
\[
\Gamma(T) \= \{ H \in \M_{n}^d \ : \ \| H \| \leq 1, \ {\rm and\ } \Re \left[ \d(\tp)^* \nabla \d(\tp) [H] \right] < 0 \} .
\]
The transverse inward set of $T$ is   the subset of $\Gamma(T)$ defined by
\[
\Delta(T) \= \{ H \in \M_{n}^d \ : \   \| H \| \leq 1, \ {\rm and\ } \d(\tp)^* \nabla \d(\tp) [H]  < 0 \} .
\]
\end{defin}

We have the following elementary result.
\begin{lemma}
\label{lemg1}
Let $H \in \Gamma(T)$. Then there exists $ \varepsilon > 0$ such that
\[
\tp + tH \in \gdel \ \forall \
0 < t < \varepsilon.
\]
Moreover, $\tp+tH$ approaches $\tp$ non-tangentially as $t \downarrow 0$.
\end{lemma}
\bp
Let $V = \d(T)$, an isometry since $T \in \idel$. Then
$$ \Re [ V^*  \nabla \d(\tp) [H]] \  \leq \ - \beta I , 
$$
 so
\se\att
\begin{eqnarray}
I - \d(\tp + tH)^* \d(\tp + tH) 
&\=&
- 2 \Re \left[ t V^*  \nabla \d(\tp) [H] \right] + O (t^2)
\label{eqg2}
 \\
 \nonumber
&\geq & 2 \beta t I + O(t^2).
\end{eqnarray}
This yields the first assertion, and the second 
follows from this and the mean value theorem, which implies that
\[
\|  \d( \tp + t H) - \d(\tp) \| \= O(t).
\]
\ep

For the rest of this section, we shall make the following assumption:

(A1) The set $\Delta(T)$ is non-empty, so there exists $K \in \mnd$ and $\beta > 0$ so that
\be
\label{eqf21}
\d(\tp)^* \nabla \d(\tp) [K] \ \leq \ - \beta I .
\ee

\begin{lemma}
\label{lemf2}
 Let $T \in \idel$, and $u$ be a $\d$ nc model for $\phi$.
Suppose  that $T$ is
  a B-point for $\phi$, and that $K \in \Delta(T)$ satisfies \eqref{eqf21}.
Let $Z_j = T + t_j K$, where $0 < t_j < 1$ and $t_j \to 0$.
If $u(Z_j)$ converges weakly to $v$, then $u(Z_j)$ converges in norm to $v$.
\end{lemma}
\bp
We have
\se
\att
\begin{eqnarray}
\label{eqf33}
I - W^* \phi(Z) &\=& v^* \left( \tensor{I}{I - \d(T)^* \d(Z) } \right) u(Z) \\
I -\phi(Z)^* \phi(Z) &\=& u(Z)^* \left( \tensor{I}{I - \d(Z)^* \d(Z) } \right) u(Z) .
\att
\label{eqf44}
\end{eqnarray}
Let $V = \d(T)$ and
 $X =  \nabla \delta(T)[K]$.
Let $Z = T + tK$, so $\d(Z) = V + tX + O(t^2) $, and recall that $V^* V = I$ since
$T \in \idel$. So the lower parts of the terms in parentheses on the right-hand sides of
\eqref{eqf33} and \eqref{eqf44} are, respectively, $-tV^* X + O(t^2)$ and $-2tV^* X + O(t^2)$.
From \eqref{eqf44} we get
\[
( u(Z)^* - v^*)  \tensor{I}{-2t V^* X} u(Z)  \=  I -\phi(Z)^* \phi(Z) - v^*\tensor{I}{-2t V^* X}u(Z) \ + O(t^2) .
\]
Take the real part,  and subtract and add twice the real part of \eqref{eqf33} to get
\beq
\Re\left[  ( u(Z)^* - v^*) \tensor{I}{-2t V^* X}  u(Z) \right]  &\=&  I -\phi(Z)^* \phi(Z) - 2 \Re(I - W^* \phi(Z)) 
+ O(t^2) \\
&=& - (I - \phi(Z)^* W)(I - W^* \phi(Z))  + O(t^2). 
\eeq
Therefore
\be
\label{eqf5}
\| \Re\left[  ( u(Z)^* - v^*) \tensor{I}{ V^* X}  u(Z) \right] \|
 \ \leq \
\frac{ \| \phi(Z) - W \|^2}{2t } + O(t).
\ee
As $Z_j   \to T$ within a non-tangential approach region,  by Theorem~\ref{thmju}
and Proposition~\ref{propf1}
\blue
$(iii) \Rightarrow (i)$, 
\black
we get 
 some constant $M$ so that
\[
\| \phi(Z_j) - W \|^2 \ \leq \
{M} \| \d(T) - \d(Z_j) \|^2 \= M \| X\|^2 t^2 + 
\blue
O(t^3) .
\black
\]
Therefore the right-hand side of \eqref{eqf5} is $O(t)$, and we conclude
that, since $V^* X \leq  - \beta I$, 
\beq
\lim_{j \to \i} 
\|   u(Z_j) -v \|^2
&\leq & 
\frac{1}{\beta}
\lim_{j \to \i} 
 ( u(Z_j)^* - v^*) \tensor{I}{ V^* X} ( u(Z_j) -v) \\
&\=& 
\frac{1}{\beta}
\lim_{j \to \i} 
\Re\left[  ( u(Z_j)^* - v^*) \tensor{I}{ V^* X}  u(Z_j)\right]  \\
 &\=& 0,
\eeq
so
\[
\lim_{j \to \i} \| v-u(Z_j)\|^2 \= 0 ,
\]
as desired.
\ep

\begin{lemma}
\label{lemf3}
Under the assumptions of Lemma~\ref{lemf2}, 
 there exists a unique $u_T$
such that 
\be
\label{eqf3}
u_T \perp {\rm ker}\left[I - \tensor{D}{I} \tensor{I}{\d(T)} \right]
\ee
and 
\be
\label{eqf4}
\left[I - \tensor{D}{I} \tensor{I}{\d(T)} \right] u_T \=
\tensor{C}{I}.
\ee
\end{lemma}
\bp
By Corollary \ref{corf2}, as $t$ decreases to $0$, the vectors $u(T + t K)$ stay bounded;
so there is some sequence $t_j$ so that $u(T + t_j K)$ converges weakly,
Choose a sequence $r_j$ \blue
increasing to $1$
\black
 so that $u(r_j T)$ converges weakly, 
and hence, by Lemma~\ref{lemf2},
 also  in norm, to a vector $v$. Writing $Z_j = T + t_j K$,
\[
\left[I - \tensor{D}{I} \tensor{I}{\d(Z_j)} \right] u(Z_j) \=
\tensor{C}{I},
\]
so taking the limit we get
\[
\left[I - \tensor{D}{I} \tensor{I}{\d(T)} \right] v \=
\tensor{C}{I}.
\]
Since $\tensor{C}{I}$ is in the range of $\left[I - \tensor{D}{I} \tensor{I}{\d(T)} \right]$,
there exists a 
 a unique vector $u_T$ satisfying \eqref{eqf4} and \eqref{eqf3}.
\ep

We can now give a characterization of B-points, for homogeneous $\delta$'s, in terms of realizations.
\bt
\label{thmf2}
Let $\phi \in \sgd$, let $u$ be a $\d$ nc model for $\phi$, and let $\begin{pmatrix}A&B\\C&D \end{pmatrix}$
be a realization as in (\ref{eqf9} -- \ref{eqf10}). Let $T \in \idel$, and assume that 
(A1) holds.
Then $T$ is a B-point for $\phi$ if and only if
\[
\tensor{C}{I} \ \in \ {\rm Ran}\left[I - \tensor{D}{I} \tensor{I}{\d(T)} \right].
\]
\et
\bp
If $T$ is a B-point, then the inclusion follows from Lemma~\ref{lemf3}.
Conversely, suppose
\[
\left[I - \tensor{D}{I} \tensor{I}{\d(T)} \right] v \= \tensor{C}{I}
\]
for some vector $v$. By \eqref{eqf1}, for any $Z$ in $\gdel$ we have
\beq
u(Z) &\=& \left[I - \tensor{D}{I} \tensor{I}{\d(Z)} \right]^{-1} \tensor{C}{I} \\
&=&
\left[I - \tensor{D}{I} \tensor{I}{\d(Z)} \right]^{-1}
\left[I - \tensor{D}{I} \tensor{I}{\d(T)} \right] v \\
&=& v 
\blue
-
\black
 \left[I - \tensor{D}{I} \tensor{I}{\d(Z)} \right]^{-1}
\left[ \tensor{D}{I} \tensor{I}{\d(T) - \d(Z)}\right] \  v.
\eeq
Then 
\se\att
\begin{eqnarray}
\nonumber
\| u(Z) \| & \ \leq \  &
\blue
\| v\| \left( 1+  \frac{\| \d(Z) - \d(T)\|}{1 - \| \d(Z)\|} \right).
\black\\
&\leq &
\| v\| \left( 1 + \frac{2 \| \d(Z) - \d(T)\|}{1 - \| \d(Z)\|^2} \right).
\label{eqf11}
\end{eqnarray}
Now let $Z$ approach $T$ non-tangentially (such as along $T+tK$), and the right-hand side of
\eqref{eqf11} stays bounded; therefore $T$ is a B-point by Corollary \ref{corf2}.
\ep

If $\delta$ is  homogeneous, then the norm of the vector $u_T$
\blue
 introduced in Lemma~\ref{lemf3}
\black
 is $\lim_{r \uparrow 1} \| u(rT) \|$.
We shall only prove it when $\delta$ is homogeneous of order $1$, though the argument
can be modified for any positive homogeneity.
\bprop
\label{propf4}
Assume that  $\delta (r Z) = r \delta(Z)$.
Let $T \in \idel$ be a B-point for $\phi$. Then $u_T$ satisfies
\se
\att
\begin{eqnarray}
\label{eqf31}
\| u_T\|^2 &\=& \lim_{r \uparrow 1} \| u(rT) \|^2 \\
\att
\label{eqf32}
&=&
\liminf_{\gdel \ni Z \to T} \frac{\| I - \phi(Z)^* \phi(Z)\|}{1 - \| \d(Z) \|^2} 
\end{eqnarray}
\eprop
\bp
Let $r_j$ be a sequence increasing to $1$ so that $u(r_j T)$ converge weakly, and hence by Lemma \ref{lemf2}
in norm, to some vector $v$.
By continuity,
\[
\left[I - \tensor{D}{I} \tensor{I}{\d(T)} \right] v \= \tensor{C}{I}, 
\]
so
\beq
u (r_j T) &\= &
\left[I - \tensor{D}{I} \tensor{I}{\d(r_j T)} \right]^{-1}\tensor{C}{I}
\\
&=& 
\left[I - \tensor{D}{I} \tensor{I}{ r_j \d(T) } \right]^{-1}
\left[I - \tensor{D}{I} \tensor{I}{\d(T)} \right] v \\
&=& 
\left[I - \tensor{D}{I} \tensor{I}{\d(T)} \right]
\left[I - \tensor{D}{I} \tensor{I}{ r_j \d(T) } \right]^{-1}
 v.
\eeq
As
$
\tensor{D}{I} \tensor{I}{\d(T)}
$
is a contraction, we have that
\[
{\rm ker}\left[ I - \tensor{D}{I} \tensor{I}{\d(T)}\right] \ \perp \
{\rm ran}\left[ I - \tensor{D}{I} \tensor{I}{\d(T)}\right].
\]
So each $u(r_j T)$ is perpendicular to ${\rm ker}\left[ I - \tensor{D}{I} \tensor{I}{\d(T)}\right]$
and hence $v$ is also. Therefore $v$ is the vector $u_T$ from 
Lemma~\ref{lemf3}, and \eqref{eqf31} holds.

To prove \eqref{eqf32}, we need to show that
\be
\label{eqf21x}
 \liminf_{\gdel \ni Z \to T} \frac{\| I - \phi(Z)^* \phi(Z)\|}{1 - \| \d(Z) \|^2} 
 \=
 \lim_{r \uparrow 1} 
 \frac{\| I - \phi(rT)^* \phi(rT)\|}{1 - \| \d(rT) \|^2} 
\ee
Let $\alpha$ denote the left-hand side of \eqref{eqf21x}. By Theorem~\ref{thmju}, we have
\beq
\frac{\| \phi(rT) - W \|^2}{\| I - \phi(rT)^* \phi(rT) \|}
&\ \leq \  &
\alpha 
\left( \frac{\| I - \d(T)^* \d(rT) \|^2} { 1 - \| \d(rT)\|^2 }\right)
\\
&=& \alpha \frac{1 - r}{1+r}.
\eeq
As in \eqref{eqbf2}, we have $\|  I - \phi(rT)^* \phi(rT) \| \leq 2 \| W - \phi(rT) \|$.
So we have
\beq
 \liminf_{r \uparrow 1} 
 \frac{\| I - \phi(rT)^* \phi(rT)\|^2}{(1 - \| \d(rT) \|^2)^2} 
 &\leq&
  \liminf_{r \uparrow 1} 
\frac{4 \|  W - \phi(rT) \|^2}{(1-r^{2})^2} \\
&\leq &
 \liminf_{r \uparrow 1}
 \frac{4}{(1-r^{2})^2} \alpha \frac{1 - r}{1+r} \| I - \phi(rT)^* \phi(rT) \|
\eeq
Dividing by 
\[
\frac{\| I - \phi(rT)^* \phi(rT) \|}{1-r^{2}}
\]
we get
\beq
 \liminf_{r \uparrow 1} 
 \frac{\| I - \phi(rT)^* \phi(rT)\|}{1 - \| \d(rT) \|^2} 
& \ \leq \ &
  \liminf_{r \uparrow 1}
  \frac{4\alpha}{(1+r)^2} \\
  &=& \alpha.
  \eeq
So we have proved that \eqref{eqf21x} holds.
\ep

\section{Derivatives at B-points}
\label{secg}

Let $T$ be a B-point of $\phi$ in $\idel$, and $W = \phi(T)$. We will keep these
fixed for the remainder of the section.
Let us make the following
assumption:

(A) The complex span of  $\Delta(T)$ is all of $\mnd$.

 This condition ensures that we have a full set of transverse directions pointing into $G_\delta.$.
 Assumption (A) is equivalent to the following two conditions holding:

(A1) The set $\Delta(T)$ is non-empty. 

(A2) The complex span of

\[
\Sigma(T) \ := \ \{ H \in \mnd \ : \ \d(\tp)^* \nabla \d(\tp) [H] \ \text{ is\ self-adjoint} \}
\]
is all of $\mnd$.

 Let $H \in \Gamma(T)$; we want to  show that
\[
\lim_{t \downarrow 0} \frac 1t [ \phi(\tp + tH) - W] 
\]
exists and is holomorphic in $H$.

First, let us sharpen Lemma \ref{lemg1}.
\begin{lemma}
\label{lemg2}
Let $ \beta > 0$. Then there exists $ \varepsilon > 0$ such that,
if
\be
\label{eqg3}
 \Re \left[ \d(\tp)^* \nabla \d(\tp) [H] \right] \ \leq \ -\beta I 
\ee
then 
\[
\tp + tH \in \gdel \ \forall \
0 < t < \varepsilon .
\]
\end{lemma}
\bp
This follows from \eqref{eqg2}, and the observation that the error term can be bounded by
some absolute constant (which depends on $\delta$ and its derivatives in a neighborhood of $T$) times $t^2$.
\ep

Let 
\[
U \= 
\{ (z,H) \in \C \times \mnd : \ H \in \Gamma(T),\ zH  \in \Gamma(T), \ T + tzH \in \gdel \ \forall \ 0 < t \leq 1 \}.
\]

Consider the set of functions 
\[
\{ \xi_t(z, H) = \frac 1t [ \phi(\tp + t z H) - W]  \ : \ 0 < t < 1 \} .
\]
These functions are all defined on $U$, and are locally bounded by Corollary~\ref{corf2}, and are holomorphic in
both $z$ and $H$.  So they form a normal family by Montel's theorem.
Let $S =( t_n)$ be  a sequence decreasing to $0$ such that
\[
\lim_{n \to \i} \xi_{t_n} (z, H) 
\]
exists; call this limit $\eta_{S}(z,H)$. We wish to show that $\eta_{S}$ does not, in fact, depend on the choice of
sequence $S$.

Let $K$ be in $\Delta(T)$. Multiplying $K$ by a small positive number if necessary, we can assume that
\[
(z,K) \in U \ \forall\ z \in \{ \Re(z) > 0 \ \text{and\ } |z| < 2 \} .
\]
Let $v$ be a unit vector in $\C^n$, and define the function 
\[
f(z) \= \la \phi(T+ zK) v, \phi(T) v \ra.
\]
Then $f : \D(1,1) \to \D$,
\blue
where
$\D(1,1)$ is the disk centered at $1$ of radius $1$.
\black
Moreover, $0$ is a B-point for $f$, because, for $t \in (0,1)$, 
\[
1 - |f(t)| \ \leq \
2\| I - \phi(T+ tK)^* \phi(T+ tK) \|
\]
and, letting $M$ denote $\| \nabla \delta (T) \|$,
\[
\| I - \delta(T+tK)^* \delta(T+tK) \| \ \leq \ 2Mt + O(t^2),
\]
so
\[
\lim_{t \downarrow 0} \frac{1 - |f(t)|}{t}
\ \leq \
 4M \ \lim_{t \downarrow 0}  \frac{\| I - \phi(T+ tK)^* \phi(T+ tK) \|}{1 - \| \d(T+tK)\|^2},
\]
and the right-hand side is bounded since $T$ is a B-point of $\phi$.

So we can apply the one variable Julia-Carath\'eodory Theorem~\ref{jc} 
to conclude that 
\[
\lim_{t \downarrow 0} \frac{ f(t) - f(0)}{t} 
\] 
exists --- indeed the limit exists as $0$ is approached non-tangentially from within $\D(1,1)$.
Since this holds for every unit vector $v$, by polarization we can conclude that
\[
\lim_{t \downarrow 0}
\frac{1}{t} [ \phi(T+ tK) -  W]
\]
exists, so every function $\eta_{S}$ agrees on points of the form $(t,K)$,
and, by holomorphicity, on points in $U$ of the form $(z,K)$, whenever $K \in \Delta(T)$.

Now, fix some element $K \in \Delta(T)$ such that $(2,K) \in U$.
Then, for some $\vare > 0$, if $H$ is in $\Sigma(T)$ and $\| H \| < \vare$,
then $K + H $ is in $ \Delta (T)$.
So all the $\eta_S$ agree on points in $U$ of the form
$(t, K+H)$, with $t > 0$ and $H$ in $\Sigma(T)$.
By assumption (A2), since $\eta_S$ is holomorphic in $H$, we get that
in fact $\eta_S$ is independent of the choice of $S$.

Thus we have proved:
\bt
\label{thmg1} Suppose $T \in \idel$ is a B-point of $\phi$, and assumption (A)  holds.
Then
\be
\label{eqg3x}
\eta(H) \= \lim_{t \downarrow 0} \frac{\phi(T + t H) - \phi(T)}{t}
\ee
exists for all $H \in \Gamma (T)$.
Moreover $\eta$ is a holomorphic function of $H$, homogeneous of degree 1.
\et

\section{Deducing the scalar case from the nc theorem}
\label{secc}

Knowing Theorem~\ref{thmju}, how could we deduce the classical Julia inequality?
We would need to know that any holomorphic function $\psi : \D \to \D$ could be
extended to a function in the Schur class of $\{ x \in \M^{[1]} : \| x \| < 1 \}$.
This indeed holds, by von Neumann's inequality \cite{vonN51}.

More generally, suppose $\Omega$ is a domain in $\C^d$, and
$\psi: \Omega \to \D$ is holomorphic. If we wish to deduce a Julia inequality using the results
of the previous section, first we need to find a matrix $\gamma$ of polynomials in $d$ commuting variables so that
$\Omega = \{ z \in \C^d: \| \gamma(z) \| < 1 \}$.
This of course may not be possible, though it is for  the polydisk, or any 
polynomial polyhedron, and for the ball.

We can define $G_\gamma$ to be the subset of ${\mathcal G}_\gamma$ consisting of {\em commuting}
$d$-tuples of matrices $x = (x^1, \dots, x^d)$ for which $\| \gamma(x) \| < 1$. The original function
$\psi$, which is holomorphic  and bounded on $G_\gamma \cap \M_1^d$, can be extended 
 to all of $G_\gamma$, either by approximating $\psi$ by polynomials, or using the Taylor functional calculus
 (see \cite{amy13a, agmc15a} for a discussion). Let us define $H^\i (G_\gamma)$ to be those
 holomorphic functions $\psi$ so that 
 \[
 \| \psi \|_{ H^\i (G_\gamma)} \ := \
 \sup \{ \| \psi(x) \| : x \in G_\gamma \}
 \]
 is finite. (If $\Omega$ were the ball, we would get the multiplier
algebra of the Drury-Arveson space).
By \cite{amfreePick15}, any function in $H^\i (G_\gamma)$ can be extended to a bounded
 function on ${\mathcal G}_\gamma$ of the same norm.
 So we can deduce the following corollary of Theorem~\ref{thmju}.
 
 \begin{cor}
 \label{corc1}
 Let $\Omega = G_\gamma \cap \M_1^d$ be a domain in $\C^d$, and let $\psi$ be a holomorphic
 function on $\Omega$. Assume that $\| \psi \|_{ H^\i (G_\gamma)} = 1$.
 Suppose $\tau \in \partial \Omega$ satisfies
 \[
 \liminf_{ z \to \tau} \frac{1 - |\psi(z)|^2}{1 - \| \gamma(z) \|^2} \= \alpha.
 \]
 Then there exists a complex number $\omega$ of modulus $1$ so that
 \[
\frac{| \psi(z) - \omega |^2}{1 - |\psi(z)|^2}
\ \leq \ 
\alpha 
\left( \frac{\| I - \gamma(\tau)^* \gamma(z) \|^2}{ 1- \| \gamma(z) \|^2}\right).
\]
 \end{cor}
 
 In the case of the ball, Corollary~\ref{corc1} is proved as Theorem 8.5.3 in \cite{rud76},
 though with the weaker assumption that $\psi$ is bounded by $1$ in the sup-norm, not in the
 multiplier norm of the Drury-Arveson space. 
 In \cite{wlo87}, K. Wlodarczyk obtained a version for the unit ball of any $J^*$-algebra, which includes
 the polydisk.

Theorem~\ref{thmg1} also can be applied to the scalar case. Assumptions (A1) and (A2) can be checked
in many concrete cases, such as the ball or the polydisk.

\begin{cor}
\label{corc2}
Assume $\Omega, \psi$ and $\tau$ are as in Corollary \ref{corc1}, and that
(A) holds at $\tau$. Then $\psi$ has a directional derivative in all inward directions at $\tau$,
and moreover this directional derivative is a holomorphic function of the direction.
\end{cor}
Of course, if $\psi$ were regular at $\tau$, the directional derivative would be a linear function of
the direction.

\section{Examples}

\begin{examp}
\label{exg1}
{\rm
Suppose
\[
\delta(Z) \=
\begin{pmatrix}
Z^1 && \\
& \ddots & \\
&& Z^d
\end{pmatrix}
\]
Then $T \in \idel$ if and only if each $T^r$ is an isometry.
We have
\[
\Delta(T) \= \{ \|K \| \leq 1 \ \text{and\ } T^{r*} K^r < 0, \ 1 \leq r \leq d \} ,
\]
and this is non-empty (\eg take $K = -T $).
The set $\Sigma(T)$ is the set of $H$ such that each $H^r$ is $T^r$ times a self-adjoint,
so (A2) is also satisfied.

Let $d=2$, and consider the scalar rational inner function
\be
\label{eqh2}
f(z,w) \= \frac{z+w-2wz}{2-z-w} .
\ee
This has a B-point at $(1,1)$. By  And\^o's theorem, $f$ is of norm one on $G_\delta$, so by
\cite{amfreePick15}, we can extend $f$ to a function of norm one on $\gdel$.
It is not immediately obvious how to do so. 

Claim: The function
\be
\label{eqh3}
\phi(Z^1, Z^2) \= \frac 12 (Z^1 + Z^2) + \frac 12 (Z^1 - Z^2) (2 - Z^1 - Z^2)^{-1} (Z^1 - Z^2)
\ee
is in $\sgd$ and agrees with $f$ on commuting variables.

Let us temporarily accept the claim.
For each $n$, the point $(I_n ,I_n)$ is a B-point, 
because
\beq
\lim_{r \uparrow 1} \frac{\| I - \phi(r I)^* \phi(r I) \|}{1-r^2} & \= & 
\lim_{r \uparrow 1} \frac{1 - r^2}{1-r^2}
 \\ &=& 1 .
\eeq
Theorem~\ref{thmju} then says that for $\phi$ as in \eqref{eqh3}, and $Z$ a pair of contractions,
\be
\label{eqh6}
\frac{\| \phi(Z) - I \|^2}{\| I - \phi(Z)^* \phi(Z) \|}
\ \leq
\ \frac{ \max_{r=1,2} \| I - Z^r \|^2}{1 - \max_{r=1,2} \| Z^r \|^2} .
\ee
If $Z^1 = Z^2$, we get equality in \eqref{eqh6}.

If we calculate $\eta(H)$ as in \eqref{eqg3x},  we get that for all $H$ with
$\Re(H^1) $ and $\Re(H^2)$ negative definite,
the derivative of $\phi$ at $(I_n, I_n) $ in the direction $H$ is
\[
\eta(H) \=  \frac 12 (H^1 + H^2)  - \frac 12 (H^1 - H^2 ) (H^1 + H^2)^{-1} (H^1 - H^2) ,
\]
which is clearly holomorphic and homogeneous of degree $1$.

Proof of claim: We shall write down a realization for $f$, and extend it to non-commutative variables.

Let
\be
\label{eqh1}
\begin{pmatrix}
A&B\\
C&D
\end{pmatrix}
\=
 \begin{pmatrix}[c|cc]
   0 & \frac{1}{\sqrt{2}} &  \frac{1}{\sqrt{2}}   \\
\hline \\
     \frac{1}{\sqrt{2}}  & \frac 12 & -  \frac 12 \\
     \frac{1}{\sqrt{2}}  & - \frac 12 &  \frac 12  \\
  \end{pmatrix} .
  \ee
Let 
\[
\delta(z,w) \= 
\begin{pmatrix}
z&0\\
0&w
\end{pmatrix},
\]
and let
\[
u(z,w) \= ( I - D \delta(z,w) )^{-1} C .
\]
Since \eqref{eqh1} is a unitary matrix, 
the function $A + B \delta(z,w) u(z,w)$ is a rational inner function on $\D^2$,
which by inspection agrees with $f$ in \eqref{eqh2}.
Now, we keep the same unitary, and using \eqref{eqf9} - \eqref{eqf10}
(where $\E$ is just $\C$, and $J=2$), we get a formula for $\phi$, which, after some algebra,
becomes
\be
\label{eqh4}
\phi(Z) \=
\frac 12 
\begin{pmatrix}
Z^1& Z^2
\end{pmatrix}
\
\begin{pmatrix}
I - \frac 12 Z^1&  \frac 12 Z^2\\
 \frac 12 Z^1&I -  \frac 12 Z^2
\end{pmatrix}^{-1}
\
\begin{pmatrix}
I\\
I
\end{pmatrix} .
\ee
Inverting the matrix on the right-hand side of \eqref{eqh4} using 
the Boltz-Banachiewicz formula does not lead to a nice formula;
but if one expands the inverse in a Neumann series, and observes that
\beq
 \begin{pmatrix}
 Z^1&  - Z^2\\
-  Z^1& Z^2
\end{pmatrix}
\
\begin{pmatrix}
I\\
I
\end{pmatrix} 
&\=&
(Z^1 - Z^2)
\begin{pmatrix}
I\\
-I
\end{pmatrix}  \\
 \begin{pmatrix}
 Z^1&  - Z^2\\
-  Z^1& Z^2
\end{pmatrix}
\
\begin{pmatrix}
A\\
-A
\end{pmatrix} 
&\=&
(Z^1 + Z^2)
\begin{pmatrix}
A\\
-A
\end{pmatrix} ,
\eeq
then for $k \geq 1$, we get
\[
\frac{1}{2^k} 
 \begin{pmatrix}
 Z^1&  - Z^2\\
-  Z^1& Z^2
\end{pmatrix}^k
\begin{pmatrix}
I\\
I
\end{pmatrix} 
\=
\frac{1}{2^k}  (Z^1 + Z^2)^{k-1} (Z^1 - Z^2) 
\begin{pmatrix}
I\\
-I
\end{pmatrix} .
\]
Then \eqref{eqh4} becomes
\[
\phi(Z) \=
 \frac 12 (Z^1 + Z^2) + \frac 14 (Z^1 - Z^2) \left[ \sum_{j=0}^\i \frac{1}{2^j} (Z^1 + Z^2)^j \right] (Z^1 - Z^2)
\]
and summing the Neumann series we get \eqref{eqh3}, as claimed.
\oec
Remark: The function
\[
\psi(Z) \= (Z^1 + Z^2 - Z^1 Z^2 - Z^2 Z^1)(2I - Z^1 - Z^2)^{-1}
\]
is another nc extension of $f$, but it is not in $\sgd$. Indeed, evaluating it on the pair of unitaries
\[
Z = \left(
\begin{pmatrix}
1 & 0\\
0 & -1
\end{pmatrix},
\
 \begin{pmatrix}
0 & 1\\
1 & 0
\end{pmatrix}
\right) ,
\]
we get
\[
\psi(Z) \= 
\begin{pmatrix}
2 & 1\\
1 &  0
\end{pmatrix},
\]
which has norm $1 + \sqrt{2}$. 
So by continuity 
\[
\| \psi \|_{\gdel} \ \geq \ 
\sup_{0 < r < 1} \| \psi(rZ) \| \ > \ 1 .
\]
}
\end{examp}
\begin{examp}
{\rm
Let $\O$ be the classical Cartan domain of symmetric $J$-by-$J$ contractive matrices in
$d = \frac{J(J+1)}{2}$ dimensions. There is an obvious embedding 
$\delta$ that takes $d$ numbers and writes them as a $J$-by-$J$ symmetric matrix, and
we can extend this map to matrices, giving $\gdel$ and the commutative version $G_\delta$.
A point is in the distinguished boundary of $\gdel$  when $ \d(T)$ is a symmetric isometry.
(A1) holds at every distinguished boundary point, and so does (A2), since
 $\Sigma(T)$ is the set of $H$ such  that $\d(H)$ can be written as the sum of $d(T)$ times a self-adjoint
 matrix and $(I - \d(T)^* \d(T))$ times anything.
So if $\phi$ is in the Schur class of $\gdel$ or $G_\delta$, we can apply both Theorem~\ref{thmju} and
Theorem \ref{thmg1}.
}\end{examp}

\bibliography{references}

\end{document}